\title{\mbox{Computing Cocycles on Simplicial Complexes\thanks{
Both authors are partially supported by the PAICYT research project
FQM-0143 from Junta de Andaluc\'{\i}a and the DGES-SEUID research project
PB97-1025-C02-02
from Education and Science Ministry (Spain).}}}
\author{Roc\'{\i}o Gonz\'{a}lez--D\'{\i}az,  Pedro Real\\
Universidad de Sevilla, Depto. de Matem\'{a}tica
Aplicada I,\\Avda. Reina Mercedes, 41012 Sevilla, Spain,\\
 e-mails: rogodi@us.es, real@us.es}
 \date{}

\documentstyle[emlines2,bezier,12pt]{article}
\newcommand{\scst}{\scriptscriptstyle}

\newtheorem{prp}{Proposition}[section]

\newtheorem{thr}[prp]{Theorem}

\newtheorem{crl}[prp]{Corollary}

\newtheorem{algorithm}[prp]{Procedure}{\bfseries}{\itshape}
\begin{document}
\maketitle
\begin{abstract}
In this note, working in the context of simplicial sets  \cite{may67},
we give a detailed study of
the complexity for computing chain level Steenrod
squares \cite{steenrod47,st52}, in terms of the number of
face operators required. This analysis is based on the
combinatorial formulation given in \cite{art1}. As an application, we give
here an algorithm for computing cup--$i$ products over integers on a
simplicial complex at chain level.
\end{abstract}
\section{Introduction}
  Cohomology operations are  tools for calculating
  $n$-cocycles on the cohomo\-logy of spaces
  (see, for example, \cite{massey,spanier}). Unfortunately,
  up to the
  present, no symbolic computational system includes {\it general}
  methods for
  finding re\-presentative $n$--cocycles on the cohomology  of spaces,
  algebras, groups, etc. Recently, several methods for finding
  $2$--cocycles representing $2$--dimensional cohomology classes of
  finite groups have been designed (see \cite{egl97,hdl93,lambe97}).
  The method established in \cite{lambe97} is
  based on the general theory presented in \cite{lambe92} and
  it seems that can be  generalized  to higher dimension without effort.

    In this paper, we describe a
  different procedure based on a combinatorial formulation given in
  \cite{art1}
   for an important class of chain
  level cohomology operations called {\em Steenrod squares}.
The formula we obtain in \cite{art1} is essentially an explicit
simplicial description of the  original formula given by Steenrod
\cite{steenrod47} for the cup--$i$ product on simplicial
complexes.
We note that a mod--$2$ explicit formulation of the Steenrod
coproduct on the chain of a simplicial set  has
also been given  in (6.2)
of Hess \cite{Hes99}, using a different method.

    We work  with simplicial sets
    \cite{may67} which are combinatorial analogs of topological
    spaces. First, our concern here is to study the ``complexity''
(in terms of number of face operators involved) of an algorithm for
computing (at chain level) the integer cup--$i$
products, using the formulation given
in \cite{art1}. Finally, as an application, we give
an algorithm for computing chain level Steenrod
squares on simplicial complexes.

We integrate here tools of Combinatorics and Computer Algebra in a
work of Algebraic Topology, opening a door to a computational development in the
search of cocycles in any degree (see \cite{casc} and \cite{gaceta}).
A treatment of some of our methods has already been presented in
\cite{aca}.

  In the literature, there   is plenty of information about
cup--$i$ products and   Steenrod squares
  (see \cite{wood98} and \cite{dieudonne89}
   for a non--exhaustive account of results).
We think that the  algorithmic technique explained here
could be substantially refined if it is suitably combined
with relevant and well--known results on these
cohomology operations and with techniques of homological
perturbation for manipulating explicit homotopy equivalences
(see  \cite{Bro,G,GL,GLS}).

We  are grateful to Prof. Julio Rubio for his helpful suggestions
for improving the algorithms showed here.
\section{Topological and Algebraic Preliminaries} \label{prel}
 The aim of this section is to give some simplicial and algebraic
 preliminaries in order to put into context
 the problem of computing $n$--cocycles (via
 cup--$i$ products and Steenrod squares). Most
of the material given in this section can be found in
 \cite{maclane75}, \cite{may67}  and \cite{spanier}.

 A {\em simplicial set} $X$ is a sequence of sets $X_0, X_1,
 \ldots$, together with {\em face operators} $\partial_i:
 X_n \rightarrow X_{n-1}$ and {\em degeneracy operators}
 $s_i: X_n \rightarrow X_{n+1}$ ($i=0,1,\ldots, n$), which
 satisfy the following simplicial identities:
\begin{eqnarray*}
\mbox{\bf (s1)} & & \partial_i \partial_j = \partial_{j-1} \partial_i \quad
\mbox{ if } i<j\enspace;\\
\mbox{\bf (s2)} & & s_i s_j= s_{j+1} s_i\quad \mbox{ if } i \leq j\enspace;\\
\mbox{\bf (s3)} & & \partial_i s_j = s_{j-1} \partial_i\quad \mbox{ if }
i<j\enspace,\\
\mbox{\bf (s4)}& &\partial_i s_j = s_j \partial_{i-1}\quad
\mbox{ if } i>j+1\enspace,\\
\mbox{\bf (s5)}& & \partial_j s_j = 1_{\scst X} =
\partial_{j+1} s_j\enspace.\end{eqnarray*}

  The elements of $X_n$ are called $n$--{\em simplices}. A simplex $x$ is
  {\em degenerate} if $x=s_i(y)$ for some simplex $y$ and
  degeneracy operator $s_i$; otherwise, $x$ is {\em non degenerate}.

Let $R$ be a ring which is commutative with unit.
Given a simplicial set $X$, let us denote $C_*(X)$ by
 the chain complex
$\{ C_{n} (X), d_n\}$, in which $C_n(X)$ is the
  free
  $R$--module generated by $X_n$ and $d_n : C_n(X) \rightarrow
  C_{n-1}(X)$ is
a $R$--module map
of degree $-1$  called  {\em differential},
defined by $d_n = \sum_{i=0}^{n} (-1)^{i} \partial_i$.

Let   $s(C_* (X))$ be the graded $R$--module generated by all the
  degenerate simplices of $X$.  In $C_*(X)$, we have that
  $d_n(s(C_{n-1}(X))) \subset s(C_{n-2}(X))$,  then \newline
  $C_{*}^{\scst N} (X)= \{ C_n (X)/s(C_{n-1} (X)), d_n\}$
is a   chain complex called {\em the normalized chain complex
  associated  to $X$}.

Since $d_n \, d_{n+1} = 0$, we can define the {\em homology} of $X$,
$H_{*}(X)$, that is the family of modules
$H_n (X)  =  \mbox{\rm Ker } d_n / \mbox{\rm Im } d_{n+1}$.

Now, the {\it cochain complex} associated to $C_*^{\scst N}(X)$,
denoted by $C^*(X;R)$, is the free $R$--module generated by all the
$R$--module maps from $C_*^{\scst N}(X)$ into $R$, together with a  map
called {\it codifferential} defined by
$(\delta^n c)(x) = c (d_{n+1} (x))$ if  $x\in C_{n+1}^{\scst
N}(X)$ and $c\in C^n(X;R)$.
We will say that $c\in C^n(X;R)$ is a $n$--{\it cocycle} if
$\delta(c)=0$, and
$c$ is a $n$--{\it coboundary} if there exists another cochain
$c'\in C^*(X;R)$ such
that $c=\delta(c')$.

 In this way, we
  define the {\it cohomology} of $X$
  with coefficients in $R$ by \newline
   $H^{*}(X) =$\,Ker\,$\delta^n/$\,Im\,$\delta^{n-1}$.
    Notice that a cocycle $c$ represents a class of cohomology.
\section{Complexity for Computing Steenrod Squares}
  First of all, let us show the explicit formula of the cup--$n$ product
$\smile_n$  on $C^{*}(X;R)$ given in \cite{art1}. The chain level
Steenrod squares $Sq^i: C^j(X;{\bf Z}_2 ) \rightarrow
C^{j+i}(X;{\bf Z}_2 )$ are defined from this operation in a very
easy way, $$Sq^i (c) = c \smile_n c, \mbox{ where } n=
j-i\enspace.$$ They satisfy that if $c$ is a $j$--cocycle, then
$Sq^i(c)$ is a $(i+j)$--cocycle.
\begin{thr}{\em \cite{art1}}\label{2}
Let $R$ be the ground ring and  $X$  a
simplicial set.  Let $c\in C^p(X;R)$, $c'\in
C^q(X;R)$ and $x\in C_{p+q-n}^{\scst N}(X)$;
if $n$ is even, then
\begin{eqnarray*}
c\smile_n c'(x)& = & \sum_{i_n=n}^{m}\quad
\sum_{i_{n-1}=n-1}^{i_n-1}\, \cdots \,
\sum_{i_0=0}^{i_1-1} \, (-1)^{A(n)+B(n,m,\bar{i})+C(n,\bar{i})+
D(n,m,\bar{i})}\\ \\
& & c(\partial_{i_0+1}\cdots \partial_{i_1-1}\partial_{i_2+1}\cdots
\cdot \partial_{i_{n-1}-1} \partial_{i_{n}+1} \cdots
\partial_{m}\,x)\\
& & \bullet c'(\partial_0 \cdots \partial_{i_0-1} \partial_{i_1+1} \cdots
\cdot \partial_{i_{n-2}-1} \partial_{i_{n-1}+1} \cdots
\partial_{i_n-1}\,x) \end{eqnarray*}
and if $n$ is odd, then
\begin{eqnarray*}
c\smile_n c'(x) &=& \sum_{i_n=n}^{m}\quad
\sum_{i_{n-1}=n-1}^{i_n-1}\, \cdots \,
\sum_{i_0=0}^{i_1-1} \ (-1)^{A(n)+B(n,m,\bar{i})+C(n,\bar{i})+
D(n,m,\bar{i})}\\\\
& & c(\partial_{i_0+1}\cdots \partial_{i_1-1}\partial_{i_2+1}\cdots
\cdot \partial_{i_{n-2}-1} \partial_{i_{n-1}+1} \cdots
\partial_{i_{n}-1}\,x) \\
& & \bullet c'(\partial_0 \cdots \partial_{i_0-1} \partial_{i_1+1} \cdots
\cdot \partial_{i_{n-1}-1} \partial_{i_n+1} \cdots
\partial_m\,x)
\end{eqnarray*}
where $m=p+q-n$, the symbol $\bullet$ is the  product in $R$,
\begin{eqnarray*}
A(n)&=&\left\{\begin{array}{l}
 1 \quad \mbox{ if } n\equiv 3,4,5,6 \mbox{ mod } 8\enspace,\\
 0 \quad \mbox{ otherwise,} \end{array}\right.\\\\
B(n,m, \bar{i})&=&\left\{\begin{array}{l}
\displaystyle\sum_{j=0}^{\lfloor\frac{n}{2}\rfloor}i_{2j} \quad \mbox{
if } n\equiv 1,2\mbox{ mod } 4\enspace,\\\\
\displaystyle\sum_{j=0}^{\lfloor\frac{n-1}{2}\rfloor}i_{2j+1}+n m \quad\mbox{ if }
n\equiv 0,3 \mbox{ mod } 4\enspace,
\end{array}\right.\\\\
C(n, \bar{i})&=&\sum_{j=1}^{\lfloor\frac{n}{2}\rfloor}
(i_{2j}+i_{2j-1})(i_{2j-1}+\cdots+i_0)
\end{eqnarray*}
and
\begin{eqnarray*}
D(n,m, \bar{i})&=&\left\{\begin{array}{l}
(m+i_n)(i_n+\cdots +i_0)\quad
\mbox{ if } n \mbox{ is odd,}\\
 0\quad\mbox{ if } n \mbox{ is even,}
\end{array}\right.
\end{eqnarray*}
being $\bar{i}= (i_0, i_1, \ldots, i_n)$.
\end{thr}

As we can see, the general organization of face operators in
these formulae is simple in the sense that we distinguish in some
way $n+1$ face operators
$\partial_{i_0},\partial_{i_1},\dots,\partial_{i_n}$; but the
signs involved follow a complicated formula. Working over
${\bf Z}_2$, this problem is eliminated.

   The aim of this section is to give an idea of the complexity
   of the algorithm for computing $n$--cocycles based  on
the previous formulation.

First of all, let us begin by giving a different description of the
cup--$n$ formula.
Let us consider an alphabet with only
two letters: $0$ and $1$. So, words in this alphabet are sequences of
letters $0$ and $1$.  We count the letters of a word
from the left to the right
and we will suppose that
the first letter on the left is in zero position.

Let $m$ and $n$ be two nonnegative integers such that
$n\leq m$. And
let $i_0,i_1,\dots,i_n\in {\bf Z}$ so that
$0\leq i_0<i_1\cdots<i_n\leq m$, then
the notation $(i_0,i_1\dots,i_n)_m$ represents the  word
with $m+1$ letters such that
there are zeros in the positions $i_0,i_1,\dots,i_n$ and ones in the rest,
that is,
$$\begin{array}{cccccccccccc}
&{\scst i_0}&&{\scst i_1}&&{\scst i_2}&&{\scst i_3}&&{\scst i_n}&\\
 1\cdots 1 &0&1\cdots 1&0&1\cdots
 1&0&1\cdots 1&0&\cdots\cdots&0&1\cdots 1\enspace.
 \end{array}$$

In the word above, by $j$--{\it block} ($1\leq j\leq n$)
we mean the block of ones
in   $i_{j-1}+1$ until $i_j-1$ positions
and zero in $i_j$ position.
The $0$--block has
ones in $0$  until $i_0-1$ positions
and zero in $i_0$ position;
and the $(n+1)$--block has
 ones in  $i_n+1$ until $m$ positions. That is,
$$\begin{array}{ccccccccc}
{\scriptstyle 0}\mbox{\scriptsize --block}&
{\scriptstyle 1}\mbox{\scriptsize --block}&
{\scriptstyle 2}\mbox{\scriptsize --block}&&
{\scriptstyle n}\mbox{\scriptsize --block}&
{\scriptstyle (n+1)}\mbox{\scriptsize --block}\\
\overbrace{\begin{array}{cc}&{\scst i_0}\\1\cdots 1&0\end{array}}&
\overbrace{\begin{array}{cc}&{\scst i_1}\\1\cdots 1&0\end{array}}&
\overbrace{\begin{array}{cc}&{\scst i_2}\\1\cdots 1&0\end{array}}&
\begin{array}{c}\\\cdots\cdots\end{array}&
\overbrace{\begin{array}{cc}&{\scst i_n}\\1\cdots1&0\end{array}}&
\overbrace{\begin{array}{c}\\1\cdots 1\end{array}}
 \end{array}$$

Eventually, the $(n+1)$--block can be the empty word.

Now, given a word  $(i_0,i_1,\dots,i_n)_m$
 we can make  a pair of words
denoted by \newline $((i_0,i_1,\dots,i_n)^+_m,(i_0,i_1\dots,i_n)^{-}_m)$,
in the following way.
If $n$ is even, then
\begin{itemize}
\item[--] the first word of the pair, denoted by
$(i_0,i_1,\dots,i_n)^{+}_m$, can be obtained from the word
$(i_0,i_1,\dots,i_n)_m$ preserving
the $j$--blocks with $j$
odd, that is,
$$\begin{array}{cccccc}
{\scriptstyle 1}\mbox{\scriptsize --bl.}&
{\scriptstyle 3}\mbox{\scriptsize --bl.}&
{\scriptstyle 5}\mbox{\scriptsize --bl.}&&
{\scriptstyle (n-1)}\mbox{\scriptsize --bl.}&
{\scriptstyle (n+1)}\mbox{\scriptsize --bl.}\\
 \overbrace{\begin{array}{c} 1\cdots1\;0 \end{array}}&
\overbrace{\begin{array}{c} 1\cdots1\;0 \end{array}}&
\overbrace{\begin{array}{c} 1\cdots1\;0 \end{array}}&
 \cdots\cdots &
\overbrace{\begin{array}{c} 1\cdots1\;0 \end{array}}&
\overbrace{\begin{array}{c} 1\cdots1 \end{array}}\enspace;
 \end{array}$$

\item[--] the second word of the pair, denoted by
$(i_0,i_1, \ldots, i_n)_m^-$ can be obtained from the word
$(i_0,i_1,\dots,i_n)_m$ preserving the
 $j$--blocks with $j$
even, that is,
$$\begin{array}{cccccccccc}
{\scriptstyle 0}\mbox{\scriptsize --bl.}&
{\scriptstyle 2}\mbox{\scriptsize --bl.}&
{\scriptstyle 4}\mbox{\scriptsize --bl.}&&
{\scriptstyle (n-2)}\mbox{\scriptsize --bl.}&
{\scriptstyle n}\mbox{\scriptsize --bl.}\\
\overbrace{\begin{array}{c}1\cdots 1\;0\end{array}}&
\overbrace{\begin{array}{c}1\cdots 1\;0\end{array}}&
\overbrace{\begin{array}{c}1\cdots 1\;0\end{array}}&
 \cdots \cdots&
\overbrace{\begin{array}{c}1\cdots 1\;0\end{array}}&
\overbrace{\begin{array}{c}1\cdots 1\;0 \end{array}}\enspace.
 \end{array}$$
\end{itemize}

   If $n$ is odd, then the procedure is analogous.

Some examples are:
\begin{itemize}
\item[--] the word $1101101$ represented by $(2,5)_6$
is associated with the pair of words:
$$((2,5)^+_6,(2,5)^-_6)=(110,1101)\enspace;$$
\item[--] the word $00110$ represented by $(0,1,4)_4$
is associated with the pair of words:
$$((0,1,4)^+_4,(0,1,4)^-_4)=(0,0110)\enspace.$$
\end{itemize}

It is easy to see that we can recover the original word
$(i_0,i_1,\ldots, i_n)_m$ from the pair $((i_0,i_1, \ldots, i_n)^+_m,
(i_0, i_1, \ldots, i_n)_m^-)$
suitably combining the $j$--blocks of both words.

For example, if we have the pair
$$(111101011, 011100)\enspace,$$
we first count the number of letters (in this case,  $m=14$), we
determine the $j$--blocks in each word of the pair
$$\begin{array}{cccc}
\qquad{\scriptstyle 0}\mbox{\scriptsize --bl.}&
{\scriptstyle 1}\mbox{\scriptsize --bl.}&
{\scriptstyle 2}\mbox{\scriptsize --bl.}\\
 \left(\quad\overbrace{\begin{array}{c} 11110 \end{array}}\right.&
\overbrace{\begin{array}{c} 10 \end{array}}&
\overbrace{\begin{array}{c} 11 \end{array}},\\
 \end{array} \;\;\;
 \begin{array}{ccc}
{\scriptstyle 0}\mbox{\scriptsize --bl.}&
{\scriptstyle 1}\mbox{\scriptsize --bl.}&
{\scriptstyle 2}\mbox{\scriptsize --bl.}\qquad\\
 \overbrace{\begin{array}{c} 0 \end{array}}&
\overbrace{\begin{array}{c} 1110 \end{array}}&
\left.\overbrace{\begin{array}{c} 0 \end{array}}\quad\right)
 \end{array}$$
and finally, we reconstruct the original word alternating the
blocks of both words
$$0\; 11110\; 1110\; 10 \; 0\; 11=(0,5,9,11,12)_{14}\enspace.$$

Identifying the letter $1$ in the position $k$ with $\partial_k$ and $0$
with the identity, the general formula for the cup--$n$ product
admits the following representation:
$$\begin{array}{ll}
c\smile_n c'(x)&\\\\
=\displaystyle\sum_{i_n=n}^{m}\,
\displaystyle\sum_{i_{n-1}=n-1}^{i_n-1}\cdots
\displaystyle\sum_{i_0=0}^{i_1-1} &(-1)^{A(n)+B(n,m, \bar{i})+
C(n,\bar{i})+D(n,m, \bar{i})}\\
&c((i_0,i_1,\dots,i_n)^+_m x) \bullet
c'((i_0,i_1,\dots,i_n)^-_m x)\enspace.\end{array}$$

And the problem of counting the number of summands in the formula
of the cup--$n$ product is equivalent to that of finding all the
possible ways to put $n+1$ zeros in $m+1$
possible places, that is,
\begin{eqnarray*}
\left(\begin{array}{c}m+1\\n+1\end{array}\right)\enspace.
\end{eqnarray*}

But, taking into account that $c$ is a $p$--cochain and $c'$ is a
$q$--cochain, then we only have to consider the summands of the
formulae having $q-n$ face operators in the first factor
and $p-n$ in the second one. Hence,
in an analogous way that in \cite{art1},
a new  combinatorial
definition of cup--$n$ product  is given in the following
theorem.
\begin{thr}\label{3}
Let $R$ be the ground ring and $X$ a
simplicial set.  If $c\in C^p(X;R)$, $c'\in
C^q(X;R)$ and $x\in C_{p+q-n}^{\scst N}(X)$, then
 \begin{eqnarray}\label{formula}\begin{array}{ll}
 c\smile_n c'(x)\\\\
= \displaystyle\sum_{i_n=S(n)}^{m}\,
\displaystyle\sum_{i_{n-1}=S(n-1)}^{i_n-1} \cdots
\displaystyle\sum_{i_1=S(1)}^{i_2-1}&(-1)^{A(n)+B(n,m,\bar{i})+C(n,\bar{i})+
D(n,m,\bar{i})}\\
&c((i_0,i_1,\dots,i_n)^+_m x)\bullet
c'((i_0,i_1,\dots,i_n)^-_m x)
\end{array}
 \end{eqnarray}
where $m=p+q-n$, $\bullet$ is the  product in $R$,
$$S(k)=i_{k+1}-i_{k-2}+\cdots+(-1)^{k+n-1}i_n+
(-1)^{k+n}\left(\lambda(n)-\left\lfloor\frac{n}{2}\right\rfloor\right)+
\left\lfloor\frac{k}{2}\right\rfloor$$
being $\lambda(n)=p$ if $n$ even and $\lambda(n)=q$ otherwise;
and $i_0=S(0)$.
\end{thr}
\noindent {\bf Proof.}

Let us start with $c\in C^p(X;R)$ and $c'\in C^q(X;R)$. If $n<p$
or $n<q$ then $c\smile_n c'$ is zero
because there is not any summand in the formula with $q-n$ face
operators in the first factor and $p-n$ face operators in the second one.
So, let us suppose that
$n\leq p$ and $n\leq q$.

If $n=0$, then $p+q-i_0=q$ and $i_0=p$, so $i_0=p$.

If $n=1$, then $i_1-1-i_0=q-1$ and $p+q-1-i_1+i_0=p-1$.
So, $i_1-i_0-q=0=q-i_1+i_0$ and hence, $i_0=i_1-q$ and $i_1\geq q$.

Let us suppose that $n$ is even (if $n$ is odd, the proof is
analogous),
then the number of face operators in the first factor of the summands
is
\begin{eqnarray}\label{factor1}
p+q-n-i_n+\cdots+i_{2k+1}-1-i_{2k}+\cdots+i_1-1-i_0\enspace,
\end{eqnarray}
and in the second one
\begin{eqnarray}\label{factor2}
i_n-1-i_{n-1}+\cdots+i_{2k}-1-i_{2k-1}+\cdots+i_2-1-i_1+i_0\enspace.
\end{eqnarray}
Since we only have to consider in the
formula for $c\smile_n c'$, the summands that the number of face
operators in the first factor is $q-n$ and $p-n$ in the second one,
that is, (\ref{factor1}) is $q-n$ and
(\ref{factor2}) is $p-n$, then
\begin{eqnarray*}
&&p+q-n-i_n+\cdots+i_{2k+1}-1-i_{2k}+\cdots+i_1-1-i_0-p+n\\
&&=i_n-1-i_{n-1}+\cdots+i_{2k}-1-i_{2k-1}+\cdots+i_2-1-i_1+i_0-q+n
\end{eqnarray*}
and hence,
\begin{eqnarray}\label{5}
i_0=i_1-i_2+i_3-\cdots-i_n+p-\frac{n}{2}\enspace.
\end{eqnarray}
Taking into account in (\ref{5}) that $i_0\geq 0$, we get
\begin{eqnarray*}
i_1\geq i_2-i_3+\cdots +i_n-p+\frac{n}{2}\enspace.
\end{eqnarray*}
Using $i_0\leq i_1-1$ in (\ref{5}), we have
\begin{eqnarray*}
i_2\geq i_3-i_4+\cdots +i_{n-1}-i_n+p-\frac{n}{2}+1\enspace.
\end{eqnarray*}
In general, let us suppose that
\begin{eqnarray*}
i_k\geq i_{k+1}-i_{k+2}+\cdots +
(-1)^{k+n-1} i_n + (-1)^{k+n}\left(p-
\frac{n}{2}\right)+\left\lfloor\frac{k}{2}\right\rfloor\enspace,
\end{eqnarray*}
for all $1\leq k\leq \ell$, and let us prove that this
expression is true in $\ell+1$ with $\ell$ odd (if $\ell$ even,
the proof is similar).
In the case $k=\ell-1$, since $i_{\ell}-1\geq i_{\ell-1}$, we have
\begin{eqnarray*}
i_{\ell}-1\geq i_{\ell}-i_{\ell+1}+\cdots +
(-1)^{\ell+n-2}i_n + (-1)^{\ell+n-1}
\left(p-\frac{n}{2}\right)+\frac{\ell-1}{2}
\end{eqnarray*}
and simplifying, we conclude
\begin{eqnarray*}
i_{\ell+1}\geq i_{\ell+2}-i_{\ell+3}+\cdots +
(-1)^{\ell+n} i_n + (-1)^{\ell+n+1}\left(
p-\frac{n}{2}\right)+\frac{\ell+1}{2}\enspace.
\end{eqnarray*}\hfill{$\Box$}
Now, let us study the number of summands in the formula
above.
Given a $p$--cochain $c$, a $q$--cochain $c'$ and a
nonnegative integer $n$,
the problem of counting all the summands in the formula
of $c\smile_n c'$ is equivalent to that of finding
all the pairs of words
$((i_0,i_1,\dots,i_n)^+_m,(i_0,i_1,\dots,i_n)^-_m)$ such that
the first  word has $q-n$ letters $1$ and the second word has $p-n$ letters
$1$. We obtain the following result.
\begin{thr}\label{tres}
Let $R$ be the ground ring. Let $X$ be a simplicial set and $n$ a
nonnegative integer. If $c\in C^p(X;R)$ and $c'\in C^q(X;R)$, then the
number of summands taking part in the formula (\ref{formula})
for $c\smile_n c'$ is
\begin{eqnarray*}
\left(\begin{array}{c}
q-\left\lfloor\frac{n+1}{2}\right\rfloor\\\\
\left\lfloor\frac{n}{2}\right\rfloor\end{array}\right)
\left(\begin{array}{c}
p-\left\lfloor\frac{n}{2}\right\rfloor\\\\
\left\lfloor\frac{n+1}{2}\right\rfloor\end{array}\right)\enspace.
\end{eqnarray*}
\end{thr}
\noindent {\bf Proof.}

First, let us suppose that $n$ is even. Our proof starts with the
observation that the first factor of a summand of
the formula (\ref{formula}) has $q-n$ face operator if
and only if the
word $(i_0,i_1,\dots,i_n)^+_m$ associated to it has
$q-n$ letters $1$ and $\frac{n}{2}$ letters $0$. Then the number of
words
$(i_0,i_1,\dots,i_n)^+_m$ having exactly $q-n$ letters $1$ is the
number of all the possible ways to put $\frac{n}{2}$ zeros
in $q-n+\frac{n}{2}$ places,
$$\left(\begin{array}{c}q-\frac{n}{2}\\\frac{n}{2}
\end{array}\right)\enspace.$$

Analogously, the word $(i_0,i_1,\dots,i_n)^-_m$  associated to the
second factor has $p-n$ letters $1$ and $\frac{n}{2}+1$ letters
$0$. Then
the number of words $(i_0,i_1,\dots,i_n)^-_m$ having $p-n$ letters $1$
 is the
number of all the possible ways to put $\frac{n}{2}$ zeros
(the last zero can not be changed)
in $p-n+\frac{n}{2}$ places, that is,
  $$\left(\begin{array}{c}p-\frac{n}{2}\\
  \frac{n}{2}\end{array}\right)\enspace.$$

And the same reasoning applied to the case $n$ odd gives us with the result
that there are
$$\left(\begin{array}{c}q-\frac{n+1}{2}\\\frac{n-1}{2}\end{array}\right)$$
possible words $(i_0,i_1,\dots,i_n)^+_m$ with $q-n$ letters $1$, and
$$\left(\begin{array}{c}p-\frac{n-1}{2}\\\frac{n+1}{2}\end{array}\right)$$
words $(i_0,i_1,\dots,i_n)^-_m$ with $p-n$ letters one.

\hfill{$\Box$}
Let us see with several examples, the improvement of the last
formulae of the cup--$n$ product given in Theorem \ref{3}
respect to the first formulae given in Theorem \ref{2}. Let us
note $c_p$ if $c\in C^p(X; R)$.
\begin{table}\caption{Number of summands}
\begin{center}
\begin{tabular}{lll}
\hline\noalign{\smallskip}
& in the formula of & in the formula of \\
& Theorem \ref{2} & Theorem \ref{3}\\
\noalign{\smallskip}\hline\noalign{\smallskip}
$c_3\smile_2 c_4$&20&6\\
$c_6\smile_5 c_6$&28&12\\
$c_{12}\smile_4 c_{10}$&11,628&1,260\\
$c_{25}\smile_5 c_{30}$&18,009,460&621,621\\
$c_{60}\smile_5 c_{70}$&4,925,156,775&68,222,616\\
$c_6\smile_5 c_{700}$&162,699,437,009,655& 970,224\\
$c_{60}\smile_{50}
c_{60}$&225,368,761,961,739,396&33,701,394,635,724,816\\
$c_6\smile_5 c_{7000}\;$&163,331,343,055,757,216,550$\;$&97,902,024\\
\hline\end{tabular}\end{center}\end{table}

Taking into account that Steenrod squares are defined
using cup--$n$ pro\-ducts, the following
corollary holds.
\begin{crl}
Let ${\bf Z}_2$ be the ground ring.
Let $i$ be a nonnegative integer and $c\in C^j(X;
{\bf Z}_2)$, then
the number of summands taking part in the
formula of $Sq^i(c)$ is
$$\left(\begin{array}{c}\left\lfloor\frac{m}{2}\right\rfloor\\\\
\left\lfloor\frac{n}{2}\right\rfloor\end{array}\right)
\left(\begin{array}{c}\left\lfloor\frac{m+1}{2}\right\rfloor\\\\
\left\lfloor\frac{n+1}{2}\right\rfloor\end{array}\right)\enspace,$$
where $m=i+j$ and $n=j-i$.
\end{crl}
\section{Simplicial Complexes}
Now, let us study a particular simplicial set.
A (combinatorial) {\em simplicial complex}
\cite{munkres,weibel} is a collection $P$ of
nonempty finite subsets of some vertex set $V$ such that if
$\tau\subset\sigma\subset V$ and $\sigma\in P$,
then $\tau\in
P$. If the vertex set is ordered, we call $P$ an {\em ordered}
simplicial complex. To every such ordered simplicial complex we
associate a simplicial set $SS(P)$ as follows. The set $SS_n(P)$
consists of all ordered $(n+1)$--tuples $\langle v_0,v_1,\dots,v_n\rangle$
 of vertices (called $n$--{\it simplices}), possibly including repetition,
 such that the underlying
set $\{v_0,v_1,\dots,v_n\}$ is in $P$ (note that
$v_0\leq v_1\leq \cdots \leq v_n$). This set is endowed with
face and degeneracy operators  defined by:
$$\partial_i\langle v_0, \dots, v_n\rangle=
\langle v_0, \dots, v_{i-1},v_{i+1},\dots, v_n\rangle$$
and
$$s_i\langle v_0, \dots, v_n\rangle=
\langle v_0, \dots, v_i, v_i, \dots, v_n\rangle\enspace.$$

Notice that a simplex  is   degenerate if it has repeated vertices;
otherwise, the simplex is  non degenerate.

Summing up, a  simplicial complex $P$ can be considered
as a combinatorial version of a triangulated polyhedron.
The strong combinatorial structure in the first one (more
precisely, in $SS(P)$) is due to
considering the degeneracy operators.

From now on, due to the fact that we will work only with ordered
simplicial complexes, we will call them  simplicial complexes, and
in order to simplify the explanation,
 we will identify the ordered simplicial complex $P$ with the
associated simplicial set $SS(P)$.
Then if $v\in P_q$,
we will say that the {\it dimension} of $v$ is $q$.
By abuse of notation, we will say that a simplex
belongs to $P$ if it belongs to $P_{\ell}$ for some $\ell$.

Let $x$ and $y$ be two simplices of $P$.
We will note $x\leq y$ if $x$ is a projection of $y$.
It is clear that a simplicial set can be given by the set of all the
simplices  with maximal dimension; and a simplex belongs to $P$ if it is a
projection of some maximal simplex of $P$.

Let $x$ and $y$ be two simplices of a simplicial complex $P$. Let us
define two operations between simplices.
Let $\{z\in P: x\leq z\mbox{ and }y\leq z\}$,
then we define $x\cup y$ as the simplex of this set
with smallest dimension (it is easy to see that $x\cup y$ is unique).
And let $\{z\in P: z\leq x \mbox{ and }z\leq y\}$, then $x\cap y$
is the simplex of this set with highest dimension (observe that $x\cap y$
is unique, too).
On the other hand, the formulation of cup--$n$ products given in
Theorem \ref{3} on a
simplicial complex is the following.
\begin{prp}\label{adap}
Let $R$ be the ground ring and
 $P$  a simplicial complex with a finite number of vertices. If
$c\in C^p(P)$ and $c'\in C^q(P)$, then for all nonne\-gative integer $n$,
$c\smile_n c'\in C^{p+q-n}(P)$ is
defined by the following formulae. Let $m=p+q-n$ and
$x=\langle v_0, \dots, v_{m}\rangle\in C_{m}(P)$, then if $n$ is even,

\begin{eqnarray*}
&c\smile_n c'(x)&\\
&=&\displaystyle\sum_{i_{n}=S(n)}^{m}\quad
\displaystyle\sum_{i_{n-1}=S(n-1)}^{i_{n}-1}\cdots
\displaystyle\sum_{i_1=S(1)}^{i_2-1}(-1)^{A(n)+B(n,m,\bar{i})+C(n,\bar{i})+
D(n,m,\bar{i})}\\\\
&&c(\langle v_0, \dots, v_{i_0}, v_{i_1}, \dots v_{i_2},v_{i_3},
\dots., v_{i_{n-2}}, v_{i_{n-1}},\dots, v_{i_n}\rangle)\\
&&\bullet c'(\langle v_{i_0}, \dots, v_{i_1}, v_{i_2}, \dots v_{i_3},v_{i_4},
\dots., v_{i_{n-1}}, v_{i_n}, \dots v_m\rangle)\enspace;
\end{eqnarray*}
and if $n$ is odd, the formula is analogous.

In these formulas, $\bullet$ is   the
product in $R$,
$$S(k)=i_{k+1}-i_{k+2}+\cdots
+(-1)^{k+n-1}i_n+(-1)^{k+n} \left\lfloor
\frac{m+1}{2}\right\rfloor +\left\lfloor\frac{k}{2}\right\rfloor$$
for all $0\leq k\leq n$, and $i_0=S(0)$.
\end{prp}
\noindent{\bf Proof.}
Using the formula from Theorem \ref{3}, we only have to notice that
$$\begin{array}{l}
\partial_0\cdots\partial_{\ell}\langle v_0, \dots,v_m\rangle=
\langle v_{\ell+1},\dots,v_m\rangle\enspace,\\
\partial_{\ell}\cdots \partial_s
\langle v_0,\dots,v_m\rangle=\langle v_0,\dots,
v_{\ell-1},v_{s+1},\dots,v_m\rangle\enspace,\\
\partial_s\cdots\partial_m\langle v_0, \dots,v_m\rangle=
\langle v_0,\dots,v_{s-1}\rangle\enspace.
\end{array}$$
\hfill{$\Box$}
For example, the formula
\begin{eqnarray*}
&&c\smile_1 c'(\langle v_0,v_1,\dots v_m\rangle)\\\\
&&=\sum_{j=0}^{p-1}(-1)^{j+(p-1+j)q}
c(\langle v_0,\dots,v_j,v_{j+q},\dots,v_m\rangle)\bullet
c'(\langle v_j,\dots,v_{j+q}\rangle)
\end{eqnarray*}
coincides with that of Steenrod given on p. 293 of
\cite{steenrod47}, up to the sign $(-1)^{p+q}$.
\section{Algorithms}
We are interested in designing algorithms for computing cocycles
using cup--$n$ products. In order to do this, we need the
following notation.

Given a simplicial complex $P$, two nonnegative integers
$n$ and $m$, and three simplices $x,y,z$ such that
$z=x\cup y=\langle v_0,\dots,v_m\rangle$ is
a $m$--simplex and
$x\cap y=\langle v_{i_0},\dots,v_{i_n}\rangle$
is a $n$--simplex, let us define the
simplices
$$\begin{array}{l}
z^{0}=\langle v_0, \dots, v_{i_0}\rangle\enspace,\\
z^{j}=\langle v_{i_{j-1}},\dots,v_{i_j}\rangle\quad\mbox{ for }
 1\leq j\leq n\enspace,\\
z^{n+1}=\langle v_{i_n},\dots,v_m\rangle\enspace.
\end{array}$$

We have the following result.
\begin{prp}\label{jarta} Let $R$ be the ground ring.
Let $P$ be a simplicial complex,
$n$  a nonnegative integer, $c\in C^p(P)$ and
 $c'\in C^q(P)$. Let $C$ (resp. $C'$) be the set of
non degenerate simplices  of $P$ such that $c(x)\neq 0$ if and only
if $x\in C$ (resp. $c'(x)\neq 0$ if and only if $x\in C'$). Let $m=p+q-n$ and
let $z=\langle v_0,\dots,v_m\rangle$ be a simplex of $P$. Define the set
$$\begin{array}{ll}
D_z=\{&(x_r,y_s): x_r\in C,\,y_s\in
C',\,x_r\cup y_s=z,\\
&x_r\cap y_s=\langle v_{i_0},\dots,v_{i_n}\rangle\mbox{ is a
$n$--simplex with }i_0=S(0),\\
&x_r=\bigcup_{j \mbox{\scriptsize even}}z^j\}\enspace.
\end{array}$$
Then,
$$c\smile_n c'(z)=\displaystyle\sum_{(x,\,y)\in
D_z}(-1)^{A(n)+B(n,m,\bar{i})+C(n,\bar{i})+D(n,m,\bar{i})}c(x)\bullet
c'(y)$$
where $\bullet$ is the product in $R$.
\end{prp}
\noindent{\bf Proof.}

Using the formula in Proposition \ref{adap} for $c\smile_n c'$,
it is not difficult to see that a summand of the formula is not zero if
the first factor is a simplex of $C$ and the second factor is a simplex
of $C'$. Hence, the simplices $x_r\in C$
and $y_s\in C'$, are both factors of a summand if and only
if $x_r\cup y_s=z$, $x_r\cap y_s=\langle v_{i_0},\dots,v_{i_n}\rangle$
is a $n$--simplex with $i_0=S(0)$ (therefore,
the rest of the inequalities $S(k)\leq i_k\leq i_{k+1}$, $1\leq k\leq
n-1$, and $S(n)\leq i_n\leq m$  are verified) and
$x_r=\bigcup_{j \mbox{ \scriptsize even}}z^j$.

\hfill{$\Box$}
Translating  this
result to a more algorithmic language, we obtain the follo\-wing
 method in which the output is expressed as a formal sum of simplices.
\begin{algorithm}Algorithm for computing cup--$n$ products.
\begin{tabbing}
{\bf Input}: \= the ground ring $R$,\\
\>  a simplicial complex $P$,\\
\> a $p$--cochain $c$ and a $q$--cochain $c'$.\\
Construct the set $C$ of $p$--simplices so that $x\in C$ if and only if
$c(x)\neq 0$.\\
Construct the set $C'$ so that $y\in C'$ if and only if $c'(y)\neq
0$.\\
 Initially, $D:=\{\;\}$.\\
 {\bf for} \= each  $x\in C$ and $y\in C'$, {\bf do}\\
 \> $z:=x\cup y=\langle v_0,\dots,v_m\rangle$,\\
 \> {\bf if} \= $x\cap y=\langle v_{i_0},\dots,v_{i_n}\rangle$ is a
$n$--simplex with $n=p+q-m$,\\
\>  \> $i_0=S(0)$ and $x=\bigcup_{j\mbox{\scriptsize even}}z^j$  {\bf then}\\
 \> \> $D:=D\cup \{(x,y)\}$.\\
 \> {\bf endif};\\
 {\bf endfor};\\
 Let $cup:=0$.\\
 {\bf for} \= each $(x,\,y)\in D_z$ {\bf do}\\
 $cup:=cup+(-1)^{A(n)+B(n,m,\bar{i})+C(n,\bar{i})+
D(n,m,\bar{i})}c(x)\bullet c(y)\; z$.\\
 {\bf endfor};\\
{\bf  Output:} \=
a formal sum, $cup=\sum \lambda_jz_j$, such that\\
\>  if $\lambda z$ is a summand of $cup$, being
$\lambda\in R$ and $z$ a $m$--simplex, then \\
\> $c\smile_n c'(z)=\lambda$, where $n=p+q-m$,\\
\> and $c\smile_n c'(z)=0$ otherwise.\\
\end{tabbing}
\end{algorithm}
Now, in order to compute  cocycles, for example,
working in ${\bf Z}_2$, we need the
following formula given in \cite{steenrod47}:
\begin{eqnarray*}
&&\delta(c\smile_n c')\\
&&=u\smile_{n-1} v+v\smile_{n-1} u+
\delta u\smile_n v+u\smile_n \delta v\enspace.
\end{eqnarray*}
It is clear that if both $c$ and $c'$ are cocycles, then
the ``commutativity" of the cup--$(n-1)$ product will determine the
obtention of cocycles via cup--$n$ products. And, in the particular case
$c=c'$, the chain level Steenrod squares appear in a
natural way.
We will develop machinery which takes advantage of this fact in a
future work.

The following result is a simple consequence of Proposition
\ref{jarta}.
\begin{crl}\label{traduccion}
Let ${\bf Z}_2$ be the ground ring. Let
$P$ be a simplicial set and $c$ a $j$--cocycle.
Let $C$ be the set of non degenerate $j$--simplices of $P$ such that
$c(x)=1$ if and only if $x\in C$. Let $i$ be a positive integer and
$z=\langle v_0,\dots,v_m\rangle$,   a $(i+j)$--simplex of $P$. Define the set
$$ \begin{array}{ll}
D_z=\{&(x_r,x_s):\, x_r,x_s\in C,\, r<s,\,
x_r\cup x_s=z,\\
&x_r\cap x_s =\langle v_{i_0},\dots, v_{i_n}\rangle
\mbox{ is a $n$--simplex with }$n=j-i$,\\
&i_0=S(0)\mbox{ and } x_r=\bigcup_{j \mbox{ \scriptsize even}}z^j
\mbox{ or }
x_r=\bigcup_{j \mbox{ \scriptsize odd}}z^j\;\}.
\end{array}$$

If cardinal of $D_z$ is even, then $Sq^i(c)(z)=0$. Otherwise,
$Sq^i(c)(z)=1$.
\end{crl}
In this corollary, we consider only $i>0$ because it is
well--known that $Sq^0(c)=c$.
And we can observe that for knowing the
cocycle $Sq^i(c)$, it is sufficient  to evaluate it over
the  simplices $z$ such
that $D_z$ is nonempty.

Using that Steenrod squares $Sq^i(c_j)$ are  cup--$(j-i)$
products, we can adapt the last algorithm to these operations. Due
to the fact that they are cohomology operations,
the first step is to determine if
  the cochain $c$ is  a cocycle.
\begin{algorithm}Algorithm for computing chain level Steenrod squares.
\begin{tabbing}
{\bf Input:} \= a simplicial complex $P$ and
 a $j$--cochain $c$.\\
Construct the set $C=\{x_1,x_2,\dots,x_k\}$
of $j$--simplices so that \\
 $x\in C$ if and only if
$c(x)=1$.\\
  Let $O:=\{\;\}$.\\
 {\bf for} \= each $x_r,\,x_s\in C$, $r<s$, and  $x_r\cup x_s$
        is a $(j+1)$--simplex {\bf do}\\
 \>  $O:= O\cup \{(x_r,\,x_s,\,x_r\cup x_s)\}$.\\
  {\bf endfor};\\
  {\bf if} \= $O$ is empty {\bf then}\\
  \> $c$ is not a cocycle\\
  {\bf else} \= let $co:=\{\;\}$.\\
  \> {\bf for} \= each $(x_r,\,x_s,\,x_r\cup x_s)\in O$ {\bf do}\\
\>  \> {\bf if} \= there exists a pair $(y,\,z)$
                in $co$ such that $x_r\cup x_s=z$ {\bf then}\\
\> \> \>   {\bf if} \= $x_r$  is not a summand of $y$, {\bf then}\\
\> \> \> \>   $co:=\left(co\setminus \{(y,\,z)\}\right)\cup \{(y+x_r,\,z)\}$.\\
\> \> \> {\bf endif};\\
\> \> \> {\bf if} \= $x_s$ is not a summand of $y$, {\bf  then}\\
\> \> \> \>  $co:=\left(co\setminus \{(y,\,z)\}\right)\cup \{(y+x_s,\,z)\}$.\\
\> \> \> {\bf endif};\\
\> \> {\bf else} $co:=co\cup \{(x_r+x_s,\,x_r\cup x_s)\}$.\\
\> \> {\bf endif};\\
\>  {\bf endfor};\\
{\bf endif};\\
 {\bf if} \= there exists some pair in $co$ such that \\
  \> the number of summands in the first element is odd {\bf then}\\
  \>   $c$ is not a cocycle,\\
  {\bf else} \= let $S:=0$.\\
 \> {\bf for} \= each $x_r,x_s\in C$, $r<s$, {\bf do}\\
 \> $z=x_r\cup x_s=\langle v_0,\dots,v_m\rangle$,\\
 \> \> {\bf if} \= $x_r\cap x_s=\langle v_{i_0},\dots,
                      v_{i_n}\rangle$ where $n=2j-m$,\\
 \> \>          \> $i_0=S(0)$\\
 \> \>          \> and $x_r=\bigcup_{t\mbox{ \scriptsize even}}z^t$
                    or $x_r=\bigcup_{t\mbox{ \scriptsize odd}}z^t$
                        {\bf then} \\
 \> \>          \>$S:=S+z$.\\
 \> \> {\bf endif};\\
 \> {\bf endfor};\\
{\bf endif}.\\
{\bf  Output:} \= a formal sum of simplices $S$
such that \\
\> if the $m$--simplex $z$ is a
summand of $S$ then \\
\> $Sq^i(c)(z)=1$ where $i=m-j$,\\
\> and $Sq^{i}(c)(z)=0$ otherwise.
\end{tabbing}
\end{algorithm}
Note that the resulting cocycles $Sq^i(c)$ can be coboundaries or not.
This means that the last algorithm does not determine Steenrod
squares at cohomology level.

The previous procedures can be easily implemented using
any Computer Algebra package or functional programming language.

\end{document}